\documentclass[11pt]{article}
\usepackage[final]{epsfig}
\usepackage{graphics}
\usepackage{amsmath}
\usepackage{amsfonts}
\usepackage{latexsym}
\usepackage{amssymb}
\usepackage{graphicx}
\title{Quasiperiodic Motion for the Pentagram Map}
\author{Valentin Ovsienko, Richard Schwartz, and Sergei Tabachnikov}

\newtheorem{theorem}{Theorem}[section]

\newtheorem{corollary}[theorem]{Corollary}

\def\P{\mbox{\boldmath{$P$}}}%
\def\PGL{\mathrm{PGL}}%
\def\R{\mbox{\boldmath{$R$}}}%
\def\Z{\mbox{\boldmath{$Z$}}}%

\begin{document}
\date{}
\maketitle

\section{Introduction and main results}

\begin{figure}[hbtp]
\centering
\includegraphics[height=1.6in]{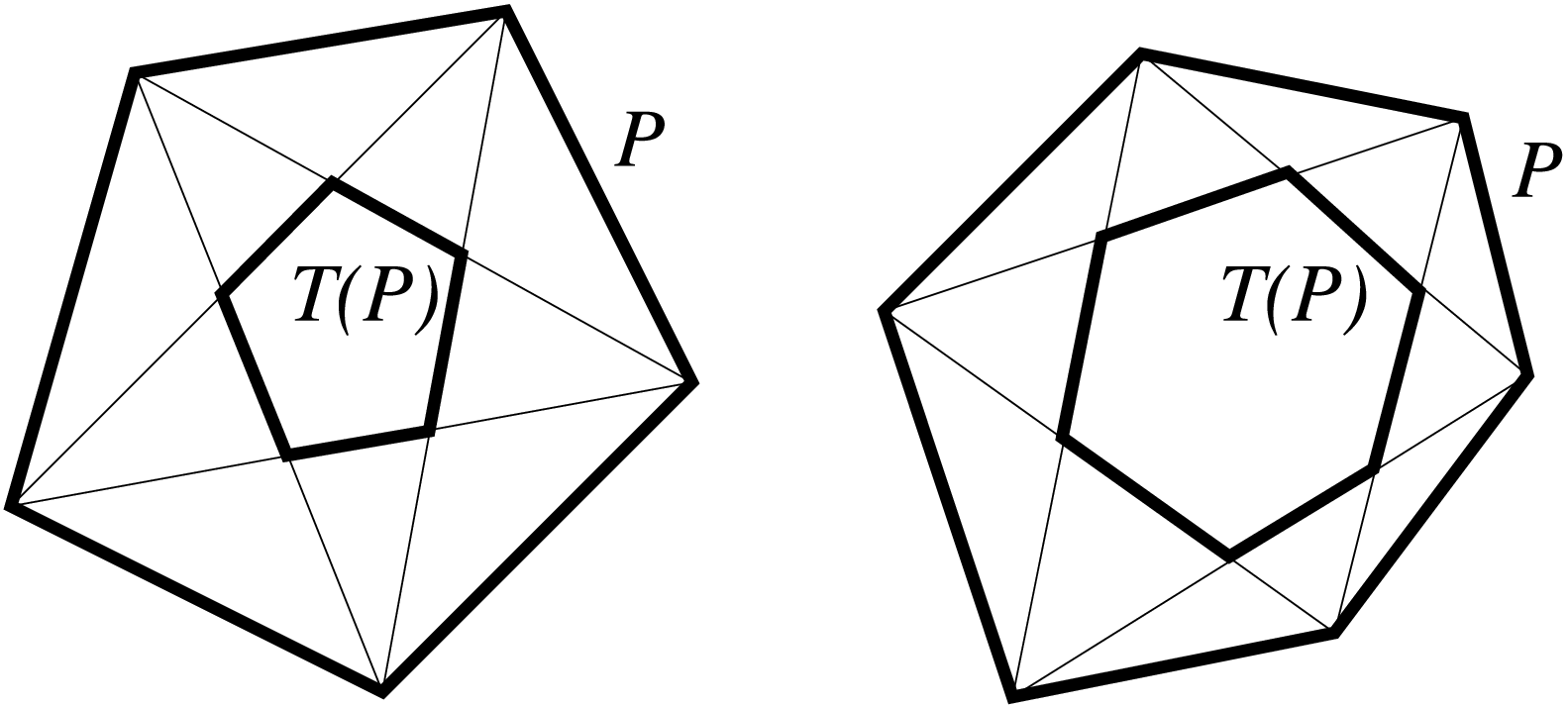}
\newline
\caption{The pentagram map defined on a pentagon and a hexagon}
\label{Fig1}
\end{figure}

The pentagram map, $T$, is a natural operation one can
perform on polygons.  See [{\bf S1\/}], [{\bf S2\/}] and [{\bf OST\/}] for
the history of this map and additional references.
Though this map can be defined for an essentially
arbitrary polygon over an essentially arbitrary
field, it is easiest to describe the map for
convex polygons contained in $\R^2$.
Given such an $n$-gon $P$, the corresponding $n$-gon $T(P)$ is the convex hull of the
intersection points of consecutive shortest diagonals of $P$.
Figure \ref{Fig1} shows two examples.

Thinking of $\R^2$ as a natural subset of the projective
plane $\R\P^2$, we observe that the pentagram map
commutes with projective transformations. 
That is, $\phi(T(P))=T(\phi(P))$, for any
$\phi\in\PGL(3,\R)$.
Let ${\cal C\/}_n$ be the space of 
\textit{convex $n$-gons modulo projective
transformations}.
The pentagram map induces a self-diffeomorphism
$T:{\cal C\/}_n\to{\cal C\/}_n$.

 $T$ is the identity map on ${\cal C\/}_5$
and an involution on ${\cal C\/}_6$, cf. [{\bf S1\/}].  
For $n \geq 7$, the
map $T$ exhibits quasi-periodic properties.
Experimentally,
the orbits of $T$ on ${\cal C\/}_n$ exhibit the kind of quasiperiodic motion
associated to a \textit{completely integrable system}.
More precisely, $T$ preserves a certain foliation of
${\cal C\/}_n$ by roughly half-dimensional tori, 
and the action of $T$ on each torus is conjugate to a
rotation.
A conjecture [{\bf S2\/}] that $T$ is completely integrable on
${\cal C\/}_n$ is still open.
 However, our recent paper [{\bf OST\/}] 
very nearly proves this result.

Rather than work directly with ${\cal C\/}_n$, we
work with a slightly larger space. 
A \textit{twisted $n$-gon} is a
map $\phi: \Z \to \R\P^2$ such that
$$
\phi(n+k)=M \circ \phi(k); \hskip 40 pt \forall k \in \Z,
$$
for some fixed element $M\in\PGL(3,\R)$ called the monodromy.
We let $v_i=\phi(i)$ and assume that
$v_{i-1},v_i,v_{i+1}$ are \textit{in general position} for all $i$. 
We denote by ${\cal P\/}_n$
the space of twisted $n$-gons modulo
projective equivalence. 
We show that the pentagram map $T:{\cal P\/}_n\to{\cal P\/}_n$
is completely integrable in the classical sense of Arnold--Liouville.
We give an explicit construction of a $T$-invariant
Poisson structure and complete list of Poisson-commuting invariants (or
integrals) for the map.  This is the algebraic
part of our theory.

The space ${\cal C\/}_n$ is naturally a subspace of
${\cal P\/}_n$, and our algebraic results say
something (but not quite enough) about the action
of the pentagram map on ${\cal C\/}_n$.  There are
still some details about how the Poisson structure
and the invariants restrict to ${\cal C\/}_n$ that
we have yet to work out.  To get a crisp geometric
result, we work with a related space, which we
describe next.

We say that a twisted $n$-polygon is {\it universally convex\/} if the
map $\phi$ is such that $\phi(\Z)\subset\R^2\subset\R\P^2$ is convex
and contained in the positive quadrant.  We also require that
the monodromy $M: \R^2 \to \R^2$ is a linear transformation
having the form
\begin{equation}
\label{hyp}
M=\left[
\begin{matrix}
a & 0 \cr 
0 & b
\end{matrix}
\right]; \hskip 40 pt  a<1<b.
\end{equation}
The image of $\phi$ looks somewhat like a ``polygonal hyperbola''.
We say that two universally convex twisted $n$-gons $\phi_1$ and $\phi_2$
are {\it equivalent\/} if there is a positive diagonal
matrix $\mu$ such that $\mu \circ \phi_1=\phi_2$.   Let
${\cal U\/}_n$ denote the space of universally convex
twisted $n$-gons modulo equivalence. It turns out that
${\cal U\/}_n$ is a pentagram-invariant and open subset of
${\cal P\/}_n$.  Here is our main
geometric result.

\begin{theorem}
\label{ac}
Almost every pont of ${\cal U\/}_n$ lies on a smooth torus that has
a $T$-invariant affine structure.  Hence, the orbit of almost every universally convex
$n$-gon undergoes quasi-periodic motion under the pentagram map.
\end{theorem}

\section{Sketch of the Proof}

In this section we will sketch the main ideas in the proof of Theorem \ref{ac}.
We refer the reader to [{\bf OST\/}] for more results and details.

\subsection{Coordinates}
\label{coord}

\begin{figure}[hbtp]
\centering
\includegraphics[height=1.8in]{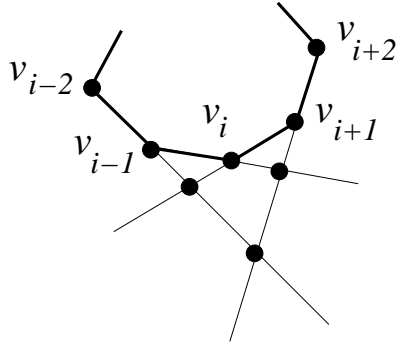}
\newline
\caption{Points involved in the definition of the invariants}
\label{Fig2}
\end{figure}

Recall that the {\it cross ratio\/} of $4$ collinear points
in $\R\P^2$ is given by
$$
[t_1,t_2,t_3,t_4]=
\frac{(t_1-t_2)\,(t_3-t_4)}{(t_1-t_3)\,(t_2-t_4)},
$$
where $t$ is (an arbitrary) affine parameter.
We use the cross ratio to construct coordinates on the
space of twisted polygons.
We associate to every vertex $v_i$ two numbers:
$$
\begin{array}{rcl}
x_i&=&
\displaystyle
\left[
v_{i-2},\,v_{i-1},\,
\left(
(v_{i-2},v_{i-1})\cap(v_i,v_{i+1})
\right),\,
\left(
(v_{i-2},v_{i-1})\cap(v_{i+1},v_{i+2})
\right)
\right]\\[10pt]
y_i&=&
\displaystyle
[\left(
(v_{i-2},v_{i-1})\cap(v_{i+1},v_{i+2})
\right),\,
\left(
(v_{i-1},v_{i})\cap(v_{i+1},v_{i+2})
\right),\,v_{i+1},\,v_{i+2}]
\end{array}
$$
called the left and right corner cross-ratios,
see Figure \ref{Fig2}.
We call our coordinates the {\it corner invariants\/}.

This construction is invariant under projective transformations,
and thus gives us coordinates on the space ${\cal P\/}_n$.
At generic points, ${\cal P\/}_n$ is locally
diffeomorphic to $\R^{2n}$.

\begin{figure}[hbtp]
\centering
\includegraphics[height=1.5in]{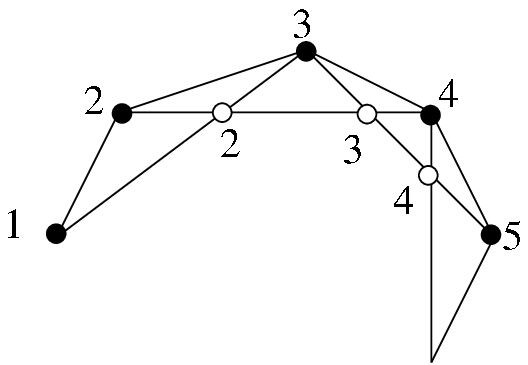}
\newline
\caption{The labelling scheme}
\label{Fig3}
\end{figure}

We will work with generic elements of ${\cal P\/}_n$,
so that all constructions are well-defined.
Let $\phi^* = T(\phi)$ be the image of
$\phi$ under the pentagram map.  We choose the
labelling scheme shown in Figure \ref{Fig3}.  The black
dots represent $\phi$ and the white ones
represent $\phi^*$.

Now we describe the pentagram map in coordinates.
\begin{equation}
\label{ExpXEq}
T^*x_i=x_i\,\frac{1-x_{i-1}\,y_{i-1}}{1-x_{i+1}\,y_{i+1}},
\qquad 
T^*y_i=y_{i+1}\,\frac{1-x_{i+2}\,y_{i+2}}{1-x_{i}\,y_{i}},
\end{equation}
Equation \ref{ExpXEq} has two immediate corollaroes.
First, there is an interesting scaling symmetry of the pentagram map.
We have a {\it rescaling operation\/} given by the
expression
\begin{equation}
R_t: \hskip 10 pt (x_1,y_1,...,x_n,y_n) \to (tx_1,t^{-1}y_1,...,tx_n,t^{-1}y_n).
\end{equation}

\begin{corollary}
The pentagram map commutes with the rescaling operation.
\end{corollary}

Second, the formula exhibits rather quickly
some invariants of the pentagram map.
For all $n$, define
\begin{equation}
\label{casimir0}
O_n=\prod_{i=1}^n x_i; \hskip 30 pt
E_n=\prod_{i=1}^n y_i
\end{equation}
When $n$ is even, define also
\begin{equation}
\label{casimir1}
O_{n/2}=\prod_{i \  {\rm even\/}} x_i+
\prod_{i \  {\rm odd\/}} x_i. \hskip 20 pt
E_{n/2}=\prod_{i \  {\rm even\/}} y_i +
\prod_{i \  {\rm odd\/}} y_i.
\end{equation}
The products in this last equation run from $1$ to $n$. 

\begin{corollary}
The functions $O_n$ and $E_n$ are
invariant under the pentagram map.
When $n$ is even, the functions $O_{n/2}$ and $E_{n/2}$ are 
also invariant under the pentagram map.
\end{corollary}

\subsection{The Monodromy Invariants}

In this section we describe the invariants of the
pentagram map.  We call them the {\it monodromy invariants\/}.
As above, let $\phi$ be a twisted $n$-gon with invariants $x_1,y_1,...$.  
Let $M$ be the monodromy of $\phi$.  We lift $M$ to an element
of $\mathrm{GL}_3(\R)$.  By slightly abusing notation, we also denote this
matrix by $M$.  The two quantities
$$
\Omega_1=\frac{{\rm trace\/}^3(M)}{{\rm det\/}(M)}; \hskip 40 pt
\Omega_2=\frac{{\rm trace\/}^3(M^{-1})}{{\rm det\/}(M^{-1})};
$$
are only dependent on the conjugacy class of $M$.

We define
$$
\widetilde \Omega_1=O_n^2E_n \Omega_1; \hskip 30 pt
\widetilde \Omega_2=O_nE_n^2 \Omega_2.
$$
In [{\bf S3\/}] (and again in [{\bf OST\/}]) it is shown that
$\widetilde \Omega_1$ and $\widetilde \Omega_2$ are 
polynomials in the corner invariants.
Since the pentagram map preserves the monodromy, and 
$O_n$ and $E_n$ are invariants, the two functions
$\widetilde \Omega_1$ and $\widetilde \Omega_2$
are also invariants.

We say that a polynomial in the corner invariants
has {\it weight\/} $k$ if
$$
R_t^*(P)=t^k P.
$$
here $R_t^*$ denotes the natural operation on polynomials
defined by the rescaling operation above.  For instance,
$O_n$ has weight $n$ and $E_n$ has weight $-n$.  In
[{\bf S3\/}] it shown that
$$
\widetilde \Omega_1=\sum_{k=1}^{[n/2]} O_k; \hskip 30 pt
\widetilde \Omega_2=\sum_{k=1}^{[n/2]} E_k
$$
where $O_k$ has weight $k$ and $E_k$ has weight $-k$.
Since the pentagram map commutes with the rescaling
operation and preserves $\widetilde \Omega_1$ and
$\widetilde \Omega_2$, it also preserves their
``weighted homogeneous parts''.   That is, the
functions $O_1,E_1,\allowbreak O_2,E_2,...$ are also invariants
of the pentagram map.  These are the monodromy
invariants.  They are all nontrivial polynomials.
In [{\bf S3\/}] it is shown
that the monodromy
invariants are algebraically independent.

The explicit formulas for the monodromy
invariants was obtained in [{\bf S3\/}].
Introduce the monomials
$$
X_i:=x_i\,y_i\,x_{i+1}.
$$
\begin{enumerate}
\item
We call two monomials $X_i$ and $X_j$ consecutive if 
$
j\in\left\{i-2,\,i-1,\,i,\,i+1,\,i+2\right\};
$
\item
we call $X_i$ and $x_j$ consecutive if
$j\in\left\{i-1,\,i,\,i+1,\,i+2\right\};$
\item
we call $x_i$ and $x_{i+1}$ consecutive.
\end{enumerate}

Let $O(X,x)$ be a monomial obtained by the product of the monomials
$X_i$ and $x_j$, i.e.,
$$
O=X_{i_1}\cdots{}X_{i_s}\,x_{j_1}\cdots{}x_{j_t}.
$$
Such a monomial is called admissible if no two of the indices 
are consecutive.
For every admissible monomial, we define the weight $|O|=s+t$
and the sign $\mathrm{sign}(O)=(-1)^t$.
One then has
$$
O_k=\sum_{|O|=k}\mathrm{sign}(O)\,O; \hskip 30 pt 
k\in\left\{1,2,\ldots,\left[\frac{n}{2}\right]\right\}.
$$
The same formula works for $E_k$, if we make all the same
definitions with $x$ and $y$ interchanged.

\subsection{The Poisson Bracket}

In [{\bf OST\/}] we introduce the Poisson bracket on $C^{\infty}({\cal P\/}_n)$.
For the coordinate functions we set
\begin{equation}
\label{PoBr}
\{x_i, x_{i\pm1}\}=\mp{}x_i\,x_{i+1},
\qquad
\{y_i,y_{i\pm1}\}=\pm{}y_i\,y_{i+1}
\end{equation}
and all other brackets vanish.
Once we have the definition on the coordinate functions,
we use linearity and the Liebniz rule to extend to all
rational functions. An easy exercise shows that
the Jacobi identity holds.

Recall the standard notions of Poisson geometry.
Two functions $f$ and $g$ are said to {\it Poisson commute\/}
if $\{f,g\}=0$.  A function $f$ is said to be a
{\it Casimir\/} (relative to the Poisson structure) if
$f$ Poisson commutes with all other functions.
The {\it corank\/} of a Poisson bracket on a 
smooth manifold is the codimension of the generic
symplectic leaves.
These symplectic leaves can be locally described as levels 
$f_i=\mathrm{const}$ of
the Casimir functions.

The main lemmas of [{\bf OST\/}] concerning our
Poisson bracket are as follows.

\begin{enumerate}
\item The Poisson bracket (\ref{PoBr})
is invariant with respect to the Pentagram map.
\item The monodromy invariants Poisson commute.
\item The invariants in Equations (\ref{casimir0}) 
and (in the even case) (\ref{casimir1}) are Casimirs.
\item The Poisson bracket has corank $2$ if $n$ if odd and corank $4$ if $n$ is even.
\end{enumerate}

We now consider the case when $n$ is odd. The even case is similar.
On the space ${\cal P\/}_n$ we have a 
generically defined and $T$-invariant
Poisson bracked that is invariant under
the pentagram map.  This bracket has co-rank $2$, and the
generic level set of the Casimir functions has dimension
$4[n/2]=2n-2$.  On the other hand, after we
exclude the two Casimirs, we have
$2[n/2]=n-1$ algebraically independent invariants that
Poisson commute with each other.  This gives us
the classical Liouville-Arnold complete integrability.

\subsection{The End of the Proof}

Now we specialize our algebraic result to the space
${\cal U\/}_n$ of universally convex twisted $n$-gons.
We check that ${\cal U\/}_n$ is an open and invariant
subset of ${\cal P\/}_n$.  The invariance is pretty clear.
The openness result derives from $3$ facts.
\begin{enumerate}
\item Local convexity is stable under perturbation.
\item The linear transformations in
Equation \ref{hyp} extend to projective transformations
whose type is stable under small perturbations.
\item A locally convex twisted polygon that has the kind of
hyperbolic monodromy given in Equation \ref{hyp} is
actually globally convex.
\end{enumerate}
As a final ingredient in our proof, we show that the
leaves of ${\cal U\/}_n$, namely the level sets of
the monodromy invariants, are compact.  We don't
need to consider all the invariants; we just
show in a direct way that the level sets of
$E_n$ and $O_n$ together are compact.

The rest of the proof is the usual application of
Sard's theorem and the definition of integrability.
We explain the main idea in the odd case.
The space ${\cal U\/}_n$ is locally diffeomorphic to
$\R^{2n}$, and foliated by leaves which generically
are smooth compact symplectic manifolds of dimension $2n-2$. A generic
point in a generic leaf lies on an $(n-1)$ dimensional smooth
compact manifold, the level set of our monodromy
invariants.  On a generic leaf, the symplectic gradients
of the monodromy functions are linearly independent at
each point of the leaf.

The $n-1$ symplectic gradients of the
monodrony invariants give a natural basis of the tangent
space at each point of our generic leaf.  This
basis is invariant under the pentagram map, and
also under the Hamiltonian flows determined by the
invariants.  This gives us a smooth compact $n-1$
manifold, admitting $n-1$ commuting flows that
preserve a natural affine structure.  From here,
we see that the leaf must be a torus.  The pentagram
map preserves the canonical basis of the torus
at each point, and hence acts as a translation.
This is the quasi-periodic motion of Theorem \ref{ac}.

\section{Pentagram map as a discrete Boussinesq equation}

Remarkably enough, the continuous
limit of the pentagram map is precisely the classical Boussinesq equation
which is one of the best known infinite-dimensional integrable systems.
This was already noticed in [{\bf S2\/}] and efficiently used in [{\bf OST\/}].
Discretization of the Boussinesq equation 
is an interesting and wel-studied subject, see [{\bf TN\/}] and references therein.
However, known versions of discrete Boussinesq equation lack geometric
interpretation.

For technical reasons we assume
throughout this section that $n\not=3\,m$.

\subsection{Difference equations and global coordinates}\label{DiffEq}

It is a powerful general idea of projective differential geometry
to represent geometrical objects in an algebraic way.
It turns out that the space of twisted $n$-gons is naturally
isomorphic to a space of difference equations.

To obtain a difference equation from a twisted polygon, lift its vertices $v_i$ to points $V_i\in\R^3$  so that $\det (V_i,V_{i+1},V_{i+2})=1$. Then
\begin{equation}
\label{diff}
V_{i+3}=a_i \,V_{i+2}+b_i\,V_{i+1} + V_i 
\end{equation}
where $a_i,b_i$ are $n$-periodic sequences of real numbers. Conversely, given two arbitrary $n$-periodic sequences $(a_i),\,(b_i)$, the difference equation (\ref{diff}) determines a projective equivalence class of a twisted polygon. This provides a global coordinate system
$(a_i,\,b_i)$ on the space of twisted $n$-gons.

\begin{corollary} 
\label{basic}
If $n$ is not divisible by 3 then the space ${\cal P\/}_n$ is isomorphic to
$\R^{2n}$.
\end{corollary}

\noindent
When $n$ is
divisible by $3$, the topology of the space is trickier.

The relation between coordinates is as follows:
$$
x_{i}=\frac{a_{i-2}}{b_{i-2}\,b_{i-1}}, 
\quad
y_{i}=-\frac{b_{i-1}}{a_{i-2}\,a_{i-1}}.
$$
The explicit formula for the pentagram map and the Poisson structure in the coordinates $(a_i,\,b_i)$ is more complicated that (\ref{ExpXEq}).
Assume $n=3m+1$ or $n=3m+2$. Then
\begin{equation}
\label{ExpXEqab}
T^*a_i=
a_{i+2}\,
\prod_{k=1}^m
\frac{1+a_{i+3k+2}\,b_{i+3k+1}}{1+a_{i-3k+2}\,b_{i-3k+1}},
\qquad 
T^*b_i=
b_{i-1}\,
\prod_{k=1}^m
\frac{1+a_{i-3k-2}\,b_{i-3k-1}}{1+a_{i+3k-2}\,b_{i+3k-1}};
\end{equation}
the Poisson bracket (\ref{PoBr}) is defined on
the coordinate functions as follows:
\begin{equation}
\label{PoBrab}
\begin{array}{rcl}
\{a_i, a_j\}&=&
\displaystyle
\sum_{k=1}^m
\left(
\delta_{i,j+3k}-\delta_{i,j-3k}
\right)a_i\,a_j,\\[16pt]
\{a_i, b_j\}&=&0,\\[6pt]
\{b_i, b_j\}&=&
\displaystyle
\sum_{k=1}^m
\left(
\delta_{i,j-3k}-\delta_{i,j+3k}
\right)b_i\,b_j.
\end{array}
\end{equation}

\noindent
The monodromy invariants also have a nice combinatorial description
in the $(a,b)$-coordinates, cf. [{\bf OST\/}].

\subsection{The continuous limit}

We understand the $n\to\infty$ continuous limit of a twisted $n$-gon  as
 a smooth parametrized curve $\gamma:\R\to\R\P^2$ with monodromy:
$$
\gamma(x+1)=M(\gamma(x)),
$$
for all $x\in\R$, where $M\in\mathrm{PGL}(3,\R)$ is fixed.
The assumption that every three consequtive
points are in general position corresponds to the assumption that
the vectors $\gamma'(x)$ and $\gamma''(x)$ are linearly independent for all $x\in\R$.
A curve $\gamma $ satisfying these conditions is usually called {\it non-degenerate}.
As in the discrete case, we consider classes of projectively equivalent curves.

The space of non-degenerate curves
is very well known in classical projective differential
geometry.
There exists a one-to-one correspondence between 
this space and the space of
linear differential operators on $\R$:
$$
A=
\left(
\frac{d}{dx}
\right)^3 
+\frac{1}{2}
\left(
u(x)\,\frac{d}{dx}+\frac{d}{dx}\,u(x)
\right)
+w(x),
$$
where $u$ and $w$ are smooth periodic functions.

\begin{figure}[hbtp]
\centering
\includegraphics[width=2in]{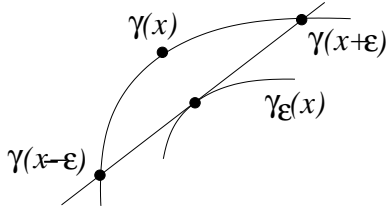}
\caption{Evolution of a non-degenerate curve}
\label{BousFig}
\end{figure}

We are looking for a continuous analog of the map $T$.
The construction is as follows.
Given a non-degenerate curve $\gamma(x)$, at each point $x$ we draw the 
chord
$\left(\gamma(x-\varepsilon), \gamma(x+ \varepsilon)\right)$ and obtain a new curve,
$\gamma_ \varepsilon(x)$, as the envelop of these chords, see Figure \ref{BousFig}.
Let $u_\varepsilon$ and $w_\varepsilon$ be the respective periodic functions.
It turns out that 
$$
u_\varepsilon=u+\varepsilon^2\widetilde{u}+(\varepsilon^3),
\quad
w_\varepsilon=w+\varepsilon^2\widetilde{w}+(\varepsilon^3),
$$
giving the curve flow: $\dot{u}=\widetilde{u},\ \dot{w}=\widetilde{w}$. We show that
$$
\displaystyle
 \dot{u}=w',\quad
 \displaystyle
\dot{w}=
 \displaystyle
 -\frac{u\,u'}{3}-\frac{u'''}{12},
$$
or
$$
\ddot{u}+\frac{\left(u^2\right)''}{6}+\frac{u^{(IV)}}{12}=0,
$$
which is nothing else but the classical Boussinesq equation.

Consider the space of functionals of the form
$$
H(u,w)=\int_{S^1}h(u,u',\ldots,w,w',\ldots)\,dx,
$$
where $h$ is a polynomial.
The \textit{first Poisson bracket} on the above space of functionals
is defined by
\begin{equation}
\label{VPBEq}
\{G,H\}=\int_{S^1}
\left(
\delta_uG\left(\delta_wH\right)'
+\delta_wG\left(\delta_uH\right)'
\right)dx,
\end{equation}
where $\delta_uH$ and $\delta_wH$ are the standard variational derivatives.
The Poisson bracket (\ref{PoBrab}) is a discrete version of
the bracket (\ref{VPBEq}).

\section*{Concluding remarks}

As it happens, this work provides more open problems than established theorems.
Let us mention here the problems that we consider most important.

\begin{enumerate}
\item
Is there there another, second $T$-invariant Poisson bracket,
compatible with the above described one?
A positive answer would allow one to apply to the pentagram map
the powerful bi-Hamiltonian techniques.
It could also help to answer the next question.
\item
Is the restriction of the pentagram map to the space
${\cal C\/}_n$ of closed polygons integrable?
\item
Perhaps the most exciting open problem is to understand
the relation of the pentagram map to cluster algebras.
It is known that the space ${\cal P\/}_n$ is a cluster manifold;
besides, our Poisson bracket has a striking similarity
with the canonical Poisson bracket on cluster manifolds [{\bf GSV\/}],
see [{\bf OST\/}] for a more detailed discussion.

\end{enumerate}
\medskip

\noindent \textbf{Acknowledgments}.
 R. S. and S. T. were partially supported by NSF grants, DMS-0604426  and DMS-0555803, respectively. V. O. and S. T. are grateful to the Research in Teams program at BIRS for its hospitality.

\section{References}

[{\bf A\/}] V. I. Arnold, Mathematical Methods of Classical Mechanics,
{\it Graduate Texts in Mathematics\/} {\bf 60\/}, Springer-Verlag, New York, 1989.
\newline
\newline
[{\bf GSV\/}] M. Gekhtman, M. Shapiro, A. Vainshtein,
{\it Cluster algebras and Poisson geometry},
Mosc. Math. J. {\bf 3} (2003), 899--934.
\newline
\newline
[{\bf OST\/}] V. Ovsienko, S. Tabachnikov, R. Schwartz, {\it The Pentagram Map: A
Discrete Integrable System\/}, submitted (2008).
\newline
\newline
[{\bf S1\/}] R. Schwartz, {\it The Pentagram Map\/}, Experimental Math. {\bf 1\/} (1992), 71--81.
\newline
\newline
[{\bf S2\/}] R. Schwartz, {\it Discrete Monodomy, Pentagrams, and the Method of
Condensation\/} J. Fixed Point Theory Appl. {\bf 3\/} (2008), 379--409.
\newline
\newline
[{\bf TN\/}]
A. Tongas, F. Nijhoff,
{\it The Boussinesq integrable system:
compatible lattice and continuum structures},
Glasg. Math. J. {\bf 47} (2005), 205--219.
\vskip 15 mm

Valentin Ovsienko: 
CNRS, Institut Camille Jordan, 
Universit\'e Lyon 1, Villeurbanne Cedex 69622, France, 
ovsienko@math.univ-lyon1.fr
\medskip

Serge Tabachnikov: Department of Mathematics, 
Pennsylvania State University,
University Park, PA 16802, USA, 
tabachni@math.psu.edu
\medskip

Richard Evan Schwartz: 
Department of Mathematics,
Brown University,
Providence, RI 02912, USA, 
res@math.brown.edu

\end{document}